\documentclass[sn-mathphys,Numbered]{sn-jnl}

\usepackage{graphicx}%
\usepackage{multirow}%
\usepackage{amsmath,amssymb,amsfonts}%
\usepackage{amsthm}%
\usepackage{mathrsfs}%
\usepackage[title]{appendix}%
\usepackage{xcolor}%
\usepackage{textcomp}%
\usepackage{manyfoot}%
\usepackage{booktabs}%
\usepackage{algorithm}%
\usepackage{algorithmicx}%
\usepackage{algpseudocode}%
\usepackage{listings}%
\usepackage{enumerate}
\usepackage{setspace}

\theoremstyle{thmstyleone}%
\newtheorem{theorem}{Theorem}[section]%
\newtheorem{lemma}{Lemma}[section]%

\theoremstyle{thmstyletwo}%

\theoremstyle{thmstylethree}%
\newtheorem{definition}{Definition}[section]%

\theoremstyle{thmstylefour}%
\newtheorem{corollary}{Corollary}[section]%

\numberwithin{equation}{section}
\begin{document}
	
	\title[Shi and Hamilton-type gradient estimates under Finsler $CD(-K,N)$ geometric flow]{Shi-type and Hamilton-type gradient estimates for a general parabolic equation under compact Finsler $CD(-K,N)$ geometric flows}
	
	
	\author[1]{\fnm{Yijie} \sur{Miao}}\email{213220672@seu.edu.cn}
	\author*[1]{\fnm{Bin} \sur{Shen}}\email{shenbin@seu.edu.cn}
	
	
	
	\affil*[1]{\orgdiv{School of Mathematics}, \orgname{Southeast University}, \orgaddress{\street{Dongnandaxue Road 2}, \city{Nanjing}, \postcode{211189}, \state{Jiangsu}, \country{China}}}


	
	\abstract{Recently, the Li-Yau-type gradient estimates for positive solutions to parabolic equations 
	\begin{equation}
		\partial_t u=\Delta u+\mathcal{R}_1u+\mathcal{R}_2u^{\alpha}+\mathcal{R}_3u(\log u)^{\beta},\notag
	\end{equation}
	under the general compact Finsler $CD(-K,N)$ geometric flow are studied. Here $\mathcal{R}_1$,$\mathcal{R}_2$,$\mathcal{R}_3$ $\in$ $ C^{1}(M,[0,T])$,  $\alpha$ and $\beta$ are both positive constants, $T$ is the maximal existence time for the flow.	
	However, compared with the Riemannian case, the curvature conditions impose stricter derivative bounds on the development term in the geometric flow, as well as on the derivative bounds of the distortion of the manifold. In this manuscript, we present Shi-type and Hamilton-type gradient estimates to demonstrate the possibility of removing such conditions.
	}

	
	
	
	\maketitle

	\section{Introduction}
	
	Gradient estimates are fundamental tools for analyzing solutions to nonlinear partial differential equations arising in geometric analysis. In the 1980s, P. Li and S.-T. Yau established the celebrated Li-Yau gradient estimate for the heat equation 
	\begin{eqnarray}\label{equ0}
		(\Delta-\partial_t)u=0,
	\end{eqnarray}
	from which they derived Harnack inequalities \cite{Yau:1986}. Building on this work, R. Hamilton introduced an elliptic-type gradient estimate, now known as the Hamilton gradient estimate \cite{Hamilton:1993}, and subsequently proved a Harnack inequality analogous to the Li-Yau result. These estimates are particularly useful for establishing monotonicity formulas. Meanwhile, W. Shi developed local derivative estimates for the Ricci flow, ensuring the existence of a solution $g_{ij}(x, t) > 0$ on the entire manifold $M$ by taking $\partial D$ to infinity on $M$ \cite{Shi:1989}. These results were later extended by B. Kotschwar \cite{Kotschwar:2007}, who established a Shi-type local gradient estimate for the heat equation on Riemannian manifolds.

    Subsequent research has generalized these ideas to broader geometric settings. In 2018, H. Dung derived a global Hamilton-type gradient estimate for a Yamabe-type parabolic equation on Riemannian manifolds with Ricci curvature bounded below \cite{Dung H T:2018}. Later, in 2021, L. Zhang investigated gradient estimates and Harnack inequalities for such equations under the Yamabe flow, further advancing the theory \cite{Zhang:2021}. More recently, H. Dung has explored various aspects of gradient estimates for positive solutions to Yamabe-type equations under different curvature conditions, along with their applications \cite{Dung H T:2025}.

    A Finsler manifold with weighted Ricci curvature bounded below naturally satisfies the curvature-dimension condition \cite{Ohta:2009}, forming what is known as a Finsler metric measure space, or a Finsler $CD(K,N)$ space when $Ric^N\geq K$. S. Ohta and K.-T. Sturm first obtained Li-Yau-type gradient estimates for positive solutions to the heat flow on compact Finsler metric measure spaces by establishing a Bochner-Weitzenb\"ock formula \cite{Ohta:2014}. This result was later generalized by Q. Xia to forward complete spaces using the Nash-Moser iteration method \cite{Xia:2023}.
	A more general nonlinear equation, encompassing both the heat equation and the Yamabe equation on a Finsler metric measure space $(M,F,\mu)$, is given by
	\begin{equation}\label{2}
		\partial_t u=\Delta u+\mathcal{R}_1u+\mathcal{R}_2u^{\alpha}+\mathcal{R}_3u(\log u)^{\beta},
	\end{equation}
	where $\Delta=\Delta_{\mu}$ is a nonlinear Laplacian associated with the measure. Recently, the second author and collaborators derived gradient estimates for the Yamabe equation and a nonlinear equation on forward complete Finsler metric measure spaces using a novel Cheng-Yau approach \cite{HLSZ:2025,SZ:2025}. Additionally, we compared the Cheng-Yau method with the Nash-Moser iteration method in studying gradient estimates for the Lichnerowicz equation \cite{MS:2024}.
	
	The Finsler geometric flow on a compact Finsler metric measure space $(M,F(t),\mu)$ is given by 
	\begin{equation}\label{1}
		\frac{\partial}{\partial t}g(x,y;t)=-2h(x,y;t),
	\end{equation}
	where $ g(x,y;t)=\frac{1}{2}\frac{\partial^2}{\partial y^i \partial y^j}F^2(x,y;t)dx^i\otimes dx^j$ is the time-dependent fundamental form of $(M,F)$, and $h(x,y;t)=h_{ij}(x,y;t)dx^i\otimes dx^j$ is a general 0-homogeneous $(0,2)$-tensor.	
	In local coordinates, the flow \eqref{1} can be expressed as 
	\[
	\frac{\partial}{\partial t}g_{ij}=-2h_{ij}.
	\]
	When $h$ is taken to be the Akbar-Zadeh Ricci tensor, i.e., $h=\mathcal{R}ic$ or $h_{ij}(x,y;t)=Ric_{ij}$, the flow \eqref{1} reduces to the Finsler Ricci flow \cite{Bao:2007}: 
	$$\frac{\partial}{\partial t}g=-2\mathcal{R}ic,$$
	where $\mathcal{R}ic(x,y)=Ric_{ij}(x,y)dx^i\otimes dx^j$, with $Ric_{ij}(x,y)=\left(\frac12\mathrm{Ric}(y)\right)_{y^iy^j}$. 
	
	Building on this framework, the second author established gradient estimates for positive solutions to the heat equation under the Finsler geometric flow in \cite{B.Shen:2025}, and further investigated the nonlinear parabolic equation \eqref{2} with X. Sun in \cite{ShenSun:2025}. The main result is as follows, which rectified the results in \cite{Azami:2023}.\\	
{\bf Theorem A.} (\cite{ShenSun:2025}, Li-Yau-type gradient estimate) {\it
Let \(\left( {M,F\left( t\right) ,m }\right)\) be an \(n\)-dimensional compact \({CD}\left( {-K,N}\right)\) solution to the Finsler geometric flow  \eqref{1}, for some \(N > n\). Suppose $u=u(x,t)$ is a positive solution to \eqref{2}, and  that there exist some positive constant \({k}_{1},{k}_{2},{k}_{3},{k}_{4},{k}_{5},{k}_{6},{k}_{7}\), such that $  -{k}_{1}g \leq  h \leq  {k}_{1}g$, \(\left| \mathcal{R}_{i}\right|  \leq  {k}_{2}\),  \(\left| {{\Delta }\mathcal{R}_{i}}\right|  \leq  {k}_{3}\), \(\left| {\nabla \mathcal{R}_{i}}\right|  \leq  {k}_{4},  i=1, 2, 3,\)  $\vert{\nabla h}\vert\leq{k}_{5}$, \(\left| {\dot{\nabla }h}\right| \leq {k}_{6}\) and ${|\nabla \tau| }  \leq  {k}_{7}$. Then, for any \(\alpha  > 1\) and \(\epsilon_{1},\epsilon_{2}>0\), we have on \(M \times  \left\lbrack  {0,T}\right\rbrack\) that
\begin{equation}\label{1.3}
	\begin{split}
		{F}^{2}\left( {\nabla \log u}\right)  - \alpha {\partial }_{t}\left( {\log u}\right)  - &\alpha \mathcal{R} _{1}- {\alpha }\mathcal{R}_{2}{u}^{a} - {\alpha }\mathcal{R}_{3}{\left( \log u\right) }^{b}\\
		&\leq \frac{N\alpha ^2}{t}+C'_1N\alpha^2+C'_2\frac{N\alpha^2}{4(\alpha-1)}+C'_3\frac{N\alpha^4}{4(\alpha-1)},
	\end{split}
\end{equation}
where the constants $C'_1,C'_2,C'_3$ are given by
\begin{align*}
	C'_{1}=&\sqrt{(1+(|a|+1)\widehat{D}+\widetilde{D}_{0}+|b|\widetilde{D}_{1})\epsilon_{1}k_{4}^{2}+(1+\widehat{D}+\widetilde{D}_{0})k_{3}+ \frac{k_{5}^{2}}{{\epsilon }_{2}}+2k_{6}^{2}+\frac{N-n+8}{4}k_{1}^{2}}\\
	&+2k_{1}+[( 2a^{2}+5|a|)\widehat{D}+2( |b(b-1)|\widetilde{D}_{2}+2|b|\widetilde{D}_{1})+|b|\widetilde{D}_{0}]k_{2},\\
	C'_{2}=&2K+\epsilon_{2}+\frac{1}{N-n}k_{7},\\
	C'_{3}=&\frac{\alpha^{2}}{\epsilon_{1}}[1+(|a|+1)\widehat{D}+\widetilde{D}_{0}+|b|\widetilde{D}_{1}].
	\end{align*}}
	Later, the second author and X. Sun also investigated Li-Yau-type gradient estimates for \eqref{2} in \cite{ShenSun:2025}, under the same curvature condition as in \cite{B.Shen:2025}. Compared to the Riemannian case (cf. \cite{Cao:2008}), the curvature condition in the Finsler setting imposes stricter derivative bounds on $h$, namely, the horizontal derivative bound $ \vert \nabla h\vert \le  L_2$ and the vertical derivative bound $\vert \dot{\nabla} h \vert \le L_3$. This raises a natural question:
	
	{\it Is it possible to relax the derivative bounds on $h$ and still derive gradient estimates for positive solutions to nonlinear parabolic equations under a Finsler geometric flow?}
	
	In this manuscript, we partially address this question by establishing Shi-type and Hamilton-type gradient estimates under the Finsler geometric flow. Specifically, we prove the following theorems:

	\begin{theorem}[Shi-type gradient estimate]\label{Shi-type} 
		Let $(M,F(t),\mu)$ be a $n$-dimensional compact $CD(-K,N)$ solution to the Finsler geometric flow \eqref{1}, for some $N > n$. Suppose $-L_1g\le h\le L_1g$, for some positive constants $L_1$.
		Let $u=u(x,t)$ be a positive solution to \eqref{2}, and assume
        $\mathcal{R}_i\le \sigma_1$,$F_{\nabla u}(\nabla^{\nabla u}\mathcal{R}_i) \le \sigma_2 $, $i=1,2,3$.
        Then, for any $\alpha>0$, $\beta >0$ and $\delta \le u\le B$, it satisfies on $M\times [0,T]$ that
		\begin{align}
			F(\nabla u)\le B^2 - \delta^2 + \mathcal{G}^{\frac{1}{2}},
		\end{align}
		where $\mathcal{G}$=$\frac{1}{4}M_0^2$+$M_1$ with 
		\begin{align*}
			M_0&=2B^2+K^{+}+L_1^{+}+ \sigma_1+\vert\alpha\vert\sigma_1\hat{D}_{0}+\sigma_1(\tilde{D}_0+\beta\tilde{D}_1),\\
			M_1&=B\sigma_2+\hat{D}_1\sigma_2+B\tilde{D}_0\sigma_2+2\sigma_1B^2+2\sigma_1\hat{D}_2+2\sigma_1B^2\tilde{D}_0,\\
			\tilde{D}_i&=\max\left\{  \vert(\log \delta )^{\beta-i}\vert,\vert(\log B)^{\beta-i}\vert\right\},\; for \;\; i=0,1,  \\
			\hat{D}_j&=\max\left\{  \delta^{\alpha-1+j},B^{\alpha-1+j}\right\}, \; for \;\; j=0,1,2. 
		\end{align*}

	\end{theorem}
	\begin{theorem}[Hamilton-type gradient estimate]\label{Hamilton-type}
	Under the same assumptions and notation as in Theorem \ref{Shi-type}, the positive solution $u=u(x,t)$ to \eqref{2} satisfies on $M\times [0,T]$ that 
    \begin{equation}
        \max\{F(\nabla \log u), F(\nabla -\log u)\}\le 2\rho\sqrt{(1+\log\frac{B}{\delta})(C_1+C_2)},
    \end{equation}
    where $\rho$ is  the reversibility given by \eqref{rho}, $\epsilon_i$ are positive constants, for $i=1,2,3$ and
    \begin{align*}
        C_1&=\frac{\sigma_2}{3}(\frac{1}{\epsilon_1}+\frac{1}{\epsilon_2}\hat{D}_0+\frac{1}{\epsilon_3}\tilde{D}_0),\\
        C_2&=\sqrt{\rho^2(K+L_1)+ \left[\vert\alpha-1\vert\sigma_1\hat{D}_0+\beta\sigma_1 \tilde{D}_1+\frac{\sigma_1}{3}(1+\hat{D}_0+\tilde{D}_0)\right]}.
    \end{align*}
	\end{theorem}
One may observe that in Theorems \ref{Shi-type} and \ref{Hamilton-type}, the condition of the bound on the horizontal derivative of the distortion, i.e., $\vert\nabla \tau\vert\leq K'$, has also been removed.

As a direct application, we also derive a classical Harnack inequality.
\begin{corollary}\label{Harnack}
Under the same assumptions and notation as in Theorem \ref{Shi-type}, it satisfies for any $x_1,x_2\in M$ and for any fixed $t_0\in(0,T)$ that
\begin{equation}
u(x_1, t_0) \le u(x_2, t_0) \cdot e^{\min\{C_3, C_4\}d_F(x_1, x_2)},
\end{equation}
where $C_3$=$\frac{B^2-\delta^2+\mathcal{G}^{\frac{1}{2}}}{\delta}$ and $C_4=2\sqrt{(1+\log\frac{B}{\delta})(C_1+C_2)},$ with $C_1,C_2$ given in Theorem \ref{Hamilton-type}.  The Finsler distance is defined by $ d_F(x_1, x_2)=\int_0^1 F(\dot{\gamma}(s))ds$, where $\gamma(s)$ is a  minimal geodesic from $x_2$ to $x_1$ in $B_p(R)$.
\end{corollary}

	\section{Preliminaries}\label{sec2}
	We introduce the background by two parts. First, we recall the Finsler metric measure space. Then, we present the related concepts on the Finsler geometric flow. Any details may find in \cite{Y.Shen:2016}.
	\subsection{Finsler metric measure space}	
	A Finsler metric measure space $(M,F,\mu)$ is a triple of a differential manifold $M$, a Finsler metric $F$ and a Borel measure $\mu$. $F$ is a positive 1-homogeneous norm $F: TM \to [ 0 , \infty )$, defined on the tangent bundle, which induces the \textit{fundamental form} $g=g_{ij}dx^i\otimes dx^j$ as $g_{i j}(x, y):=\frac{1}{2} \frac{\partial^2 F^2}{\partial y^i \partial y^j}(x, y)$. Derivative along the fibre again yields the \textit{Cartan tensor}, that is
	\begin{equation*}
		C_y(X,Y,Z):=C_{ijk}(y)X^iY^jZ^k:=\frac{1}{4}\frac{\partial^3 F^2(x,y)}{\partial y^i\partial y^j\partial y^k}X^iY^jZ^k,   
	\end{equation*}
	for any local vector fields $X,Y,Z$. There is a unique almost $g$-compatible and torsion-free connection on the pull back tangent bundle $\pi^{*} TM$, called the \textit{Chern connection}, which is defined by
	\begin{align*}
		\nabla_XY-\nabla_YX&=[X,Y]; \notag  \\
		Z(g_y(X,Y))-g_y(\nabla_ZX,Y)-g_y(X,\nabla_ZY)&=2C_y(\nabla_Zy,X,Y).
	\end{align*}
	The Chern connection
	coefficients are denoted by $\Gamma_{jk}^i(x, y)$ in the natural coordinate system, which could
	induce locally the \textit{spray coefficients} as $G^i = \frac{1}{2}\Gamma_
	{jk}^iy^jy^k$. The spray is given by
	\begin{equation*}
		G=y^i\frac{\delta}{\delta x^i}=y^i\frac{\partial}{\partial x^i}-2G^i\frac{\partial}{\partial y^i}, 
	\end{equation*} 
	in which $\frac{\delta}{\delta x^{i}}=\frac{\partial}{\partial x^{i}}-N_{i}^{j}\frac{\partial}{\partial y^{j}}$, and the nonlinear connection coefficients are locally induced
	from the spray coefficients by $N_{j}^{i}=\frac{\partial G^{i}}{\partial y^{j}}$. Customarily, we denote the horizontal Chern derivative by “ 
	$\mid$ ” and the vertical Chern
	derivative by “ ; ”. For example,
	\begin{equation*}
	    w_{i|j}=\frac{\delta}{\delta x^{j}}w_{i}-\Gamma_{ij}^{k}w_{k},\quad w_{i;j}=F\frac{\partial}{\partial y^{j}}w_{i}, 
	\end{equation*}
	for any 1-form $w = w_idx^i$ on the pull-back bundle. In this manuscript, the horizontal Chern derivative also denoted by “$\nabla$”. 
	
	The \textit{Chern Riemannian curvature} $R$ is locally defined by
	\begin{equation*}
		R_{j\;kl}^{i}=\frac{\delta\Gamma_{jl}^{i}}{\delta x^{k}}-\frac{\delta\Gamma_{jk}^{i}}{\delta x^{l}}+\Gamma_{km}^{i}\Gamma_{jl}^{m}-\Gamma_{lm}^{i}\Gamma_{jk}^{m}.
	\end{equation*}
	The \textit{flag curvature} with pole $y$ is defined by
	\begin{equation*}
		K(P,y)=K(y,u):=\frac{R_{y}(y,u,u,y)}{F^{2}(y)h_{y}(u,u)}=\frac{-R_{ijkl}(y)y^iu^jy^ku^l}{(g_{ik}(y)g_{jl}(y)-g_{il}(y)g_{jk}(y))y^iu^jy^ku^l},  
	\end{equation*}
	for any two linearly independent vectors $y$, $u$ $\in$ $T_xM \setminus \{0\}$, which span a tangent plane $\Pi_{y}=\operatorname{span}\{y,u\}$. Then the \textit{Finslerian Ricci curvature} is defined by 
	\begin{equation}
		Ric(y):=F^2(y)\sum_{i=1}^{n-1}K(y,e_i),   \notag
	\end{equation}
	where $e_1,\cdots, e_{n-1}, \frac{y}{F(y)}$ form an orthonormal basis of $T_xM$ with respect to $g_y$.

	A Finsler structure is defined to be \textit{reversible} if $F(x,V)=\bar{F}(x,V)$ for any point $x$ and any vector field $V$, where  $\bar{F}(x,V):=F(x,-V)$ is called the \textit{reversed Finsler metric} of $F$. We define the \textit{reversibility} of $(M,F,\mu)$ by 
	\begin{equation}\label{rho}
		\rho:=\sup_{x\in M}\sup_{V\neq 0}\frac{F(x, V)}{\bar{F}(x, V)}.
	\end{equation}
	The forward distance from $p$ to $q$ is defined by
	\begin{equation}
		d_p(q):=d(p, q):=\inf _\gamma \int_0^1 F(\gamma(t), \dot{\gamma}(t)) dt,   \notag
	\end{equation}
	where the infimum is taken over all the $C^1$
	curves $\gamma$ : [0, 1] $\rightarrow$ $M$ such that $\gamma(0) = p$
	and $\gamma(1) = q$. A $C^2$ curve $\gamma$ is called a \emph{geodesic} if it locally satisfies the geodesic eqaution
	$$\ddot\gamma^i(t)+2G^i(\gamma(t),\dot \gamma(t))=0.$$
	
	The \textit{Legendre transformation} is an isomorphism $\mathcal{L}^*:T_x^*M\to T_xM$ maps $v^*\in T_x^*M$ to the unique element $v \in T_xM$ such that $v^*( v) = F( v) $. For a differentiable function $f:M \to \mathbb{R}$, the \textit{gradient~vector} of $u$ at $x$ is defined by $\nabla( x) : = \mathcal{L} ^* ( du( x) ) \in T_xM.$ If $du( x) = 0, $ then we have $\nabla u( x) = 0.$ If $Du( x) \neq 0, $ and we can write on $M_u:=\{x\in M:du(x)\neq0\}$ that
	\begin{equation}
		\nabla u=\sum_{i,j=1}^{n}g^{ij}(x,\nabla u)\frac{\partial u}{\partial x^{i}}\frac{\partial}{\partial x^{j}}.   \notag
	\end{equation}
	In local coordinates 
	$\{x^{i}\}_{i=1}^{n}$, expressing the volume measure by $d\mu=e^{\Phi}dx^{1}\wedge\cdots\wedge dx^{n}$, then the \textit{divergence} of a smooth vector field $V$ can be written locally as
	\begin{equation}
		\mathrm{div}_{\mu}V=\sum_{i=1}^{n}(\frac{\partial V^{i}}{\partial x^{i}}+V^{i}\frac{\partial\Phi}{\partial x^{i}}).  
	\end{equation}
	The \textit{Finsler Laplacian} of a function $f$ on $M$ could also be given by
	\begin{equation}
		\Delta f:=\mathrm{div}_{\mu}(\nabla f),   \notag
	\end{equation}
	in the distributional sense. That is, 
	\begin{equation*}
		\int_{M}\phi\Delta_{\mu}fd\mu=-\int_{M}d\phi(\nabla f)d\mu,  
	\end{equation*}
	for any $f \in W^{1,p}(M)$ and a test function $\phi \in C_0^{\infty}(M)$.
	
	The distortion of $(M, F, \mu)$ and the \textit{S-curvature} are defined by
	\begin{equation}
		\tau(x,y):=\frac12\log\det g_{ij}(x,y)-\Phi(x), \mbox{ and } S(x,y):=\frac{d}{dt}(\tau \circ \gamma)|_{t=0},
	\end{equation}
	respectively, where $\gamma = \gamma(t)$ is a forward geodesic from $x$ with the initial tangent vector
	$\dot{\gamma}(0) = y$.

	The following \textit{weighted Ricci curvature}, which is proved to be equal to the $CD(-K,N)$ condition. 	
	\begin{definition}[\cite{Ohta:2009}]
		Given a unit vector $V \in T_xM$ on a Finsler metric measure space $(M,F,\mu)$ and a positive number $k$,
		the weighted Ricci curvature is defined by		
		\begin{equation}
			Ric^k(V):=\begin{cases}Ric(x,V)+\dot{S}(x,V),&\text{if }S(x,V)=0\text{ and }k=n\text{ or if }k=\infty;\\-\infty,&\text{if }S(x,V)\neq0\text{ and if }k=n;\\Ric(x,V)+\dot{S}(x,V)-\frac{S^2(x,V)}{k-n},&\text{if }n<k<\infty,\end{cases}   \notag
		\end{equation}
		where $\dot{S}(x,V)$ is the derivative along the geodesic from $x$ in the direction of $V$ .
	\end{definition}
	Based on it, we define that
	\begin{definition}[\cite{B.Shen:2025}]
		A Finsler metric measure space satisfying $\mathrm{Ric}^{N}$ bounded from below for some $N>n$ by $-K$ with $K>0$ is called a \textit{Finsler $CD(-K,N)$ space}.
	\end{definition}
	The following Bochner-Weitzenb\"ock formula is adopted to generalize the Li-Yau gradient estimates on Finsler $CD(-K,N)$ spaces. 
	\begin{lemma}[Bochner–Weitzenb\"ock formula \cite{Ohta:2014}]
		Let $(M, F, \mu)$ be an n-dimensional Finsler metric measure space. Given $u$ $\in$ $H_{loc}^2(M) \cap C^1(M)$  with $\Delta u \in H_{loc}^1(M)$, we have 
		\begin{equation} \label{lem2.1}
			\Delta^{\nabla u}(F^2(\nabla u))=2du(\nabla^{\nabla u} (\Delta u))+2\Vert \nabla^2u\Vert^2_{HS(\nabla u)}+2 Ric^{\infty}(\nabla u).
		\end{equation}
	\end{lemma}
		
	\subsection{Finsler geometric flow}
	 In this subsection, we will give the notations utilized in the analysis of nonlinear equations under the Finsler geometric flow. Considering the Finsler geometric flow and denoting $h(y)=h_{ij}(y)y^iy^j$, \eqref{1} can be rewritten as 
	\begin{eqnarray}
		\frac{\partial}{\partial t}\log F=-H,
	\end{eqnarray}
	where $H=h(y)/F^2$ is a function defined on the sphere bundle $SM$. So the solution to the Finsler geometric flow is a time-dependent Finsler metric $F(x,y;t)$. It can be verified directly in local coordinates that 
    \begin{eqnarray}
    	\frac{\partial}{\partial t}g^{ij}=2h^{ij},
    \end{eqnarray} 
    where $(g^{ij}(x,y;t))$ is the inverse matrix of the fundamental form and $h^{ij}(x,y;t)=g^{ik}h_{kl}g^{lj}$. 
    
    Fixed a local vector field $V$ on $M$, $h(x,V;t)=h_{ij}(V)dx^i\otimes dx^j$ provides a time-dependent symmetrical bilinear form, that is,
    \begin{eqnarray}\label{a2.6}
    	h_V(W_1,W_2)=h_{ij}(V)W_1^iW_2^j,
    \end{eqnarray}
    for any two vector fields $W_1$ and $W_2$. Obviously, $h_V(V,V)=h(V)$. 
    When $h_V(\cdot,\cdot)$ acts on a (0,2)-type tensor $\Omega$, we just denote it by $h_V(\Omega)$. 
    
    For more information on the calculation and properties of Finsler geometric flow, please refer to \cite{B.Shen:2025}.

	\section{Shi-type gradient estimate}\label{sec3}
		In this section, we give a lemma about functions on a time-dependent Finsler metric measure space, establish a Shi-type gradient estimate for positive solutions to equation \eqref{2}  under the Finsler geometric flow \eqref{1}, and further derive a Harnack inequality.
	
	\begin{lemma}\label{lem1}
		Let $(M,F(t),\mu)$ be a solution to the Finsler geometric flow \eqref{1}. For any $f\in C^1(M)\cap C^1([0,T])$, it satisfies that
		\begin{equation*}
			\partial_t(F^2(\nabla f))=2h(\nabla f)+2df_{t}(\nabla f),
		\end{equation*}
		where $h(\nabla f)=h_{\nabla f}(\nabla f, \nabla f)=h_{ij}(\nabla f)f^if^j$.
	\end{lemma}
	\begin{proof}
		Direct computation gives that
		\begin{align*}
			\partial_t(F^2(\nabla f))=&\partial_t(g^{ij}(\nabla f)f_if_j)\notag\\
			=&(\partial_tg^{ij})(\nabla f)f_if_j+\frac{\partial g^{ij}}{\partial y^k}\frac{\partial f^k}{\partial t}f_if_j+2g^{ij}(\nabla f)\frac{\partial f_i}{\partial t}f_j\notag\\
			=& 2h^{ij}(\nabla f)f_if_j-2C_{\;k}^{ij}(\nabla f)\frac{\partial f^k}{\partial t}f_if_j+2g^{ij}(\nabla f)\frac{\partial f_i}{\partial t}f_j\notag\\
			=& 2h^{ij}(\nabla f)(\nabla f, \nabla f)+2df_{t}(\nabla f).
		\end{align*}
	\end{proof}
	We now turn to the positive solution $u=u(x,t)$ to the equation \eqref{2} under the Finsler geometric flow.
  
	\begin{proof}[Proof of Theorem \ref{Shi-type}]

By denoting $H=F(\nabla u)+u^2$, we will investigate the behaviour of $H$ under the quasilinear parabolic operator 
\begin{equation}\label{L}
\mathcal{L}=\Delta^{\nabla u}-\partial_t.
\end{equation}
Firstly, by applying the Bochner-Weitzenb\"ock formula, we compute $\Delta^{\nabla u} F^2(\nabla u)$ as
\begin{equation}
\Delta^{\nabla u} F^2(\nabla u)=2\Vert \nabla^2u \Vert_{HS(\nabla u)}^2+2Ric^{\infty}(\nabla u)+2d(\Delta u)(\nabla u).
\end{equation}
According to Lemma \ref{lem1}, we obtain 
\begin{equation}\label{equ:6}
\partial_tF^2(\nabla u)=2h(\nabla u)(\nabla u,\nabla u)+2du_t(\nabla u).
\end{equation}
Therefore,
\begin{align}\label{equ:3}
\mathcal{L}F^2(\nabla u)=&2\Vert \nabla^2u \Vert_{HS(\nabla u)}^2+2Ric^{\infty}(\nabla u)+2d(\Delta u)(\nabla u) \notag \\
&-2h(\nabla u)(\nabla u,\nabla u)-2du_t(\nabla u).
\end{align}
Moreover,
\begin{align} \label{equ:4}
\Delta^{\nabla u} F^2(\nabla u)&=div(\nabla^{\nabla u} (F^2(\nabla u)))\notag\\
&=2F(\nabla u)\Delta^{\nabla u}(F(\nabla u))+2F_{\nabla u}^2(\nabla^{\nabla u} F(\nabla u)).
\end{align}
Hence,
 \begin{equation}\label{equ:5}
\mathcal{L}F^2(\nabla u)=2F(\nabla u)\mathcal{L}(F(\nabla u))+2F_{\nabla u}^2(\nabla^{\nabla u} F(\nabla u)).
\end{equation}
Combining \eqref{equ:3} and \eqref{equ:5}, we have 
\begin{align}
2F(\nabla u)\mathcal{L}(F(\nabla u))+2F_{\nabla u}^2(\nabla^{\nabla u} F(\nabla u))&=2\Vert \nabla^2u \Vert_{HS(\nabla u)}^2+2Ric^{\infty}(\nabla u)+2d(\Delta u)(\nabla u) \notag \\
&\quad-2h(\nabla u)(\nabla u,\nabla u)-2du_t(\nabla u).
\end{align}
Kato's inequality on Finsler manifold could be established by direct computation in local coordinates that
\begin{align}
F_{\nabla u}^2(\nabla^{\nabla u} F(\nabla u))
&=(g^{kl}u_{j|k}u_{q|l})u^j u^q \notag \\
&\le (g^{kl}u_{j|k}u_{q|l}g^{jq})(g_{jq}u^ju^q)=\Vert\nabla^2u\Vert_{HS(\nabla u)}^2.\notag
\end{align}
Based on the assumptions $-L_1g<h<L_1g$ and the $CD(-K,N)$ condition, we derive
\begin{align}\label{equ:7}
2F(\nabla u)\mathcal{L}(F(\nabla u))&\ge 2Ric^{\infty}(\nabla u)+2d(\Delta u)(\nabla u)-2h(\nabla u)(\nabla u,\nabla u)-2du_t(\nabla u)\notag\\
&\ge -2(K+L_1)F^2(\nabla u)-2d\left(\mathcal{R}_1u+\mathcal{R}_2u^{\alpha}+\mathcal{R}_3u(\log u)^{\beta}\right)(\nabla u).
\end{align}
Next, we compute the term $d\left(\mathcal{R}_1u+\mathcal{R}_2u^{\alpha}+\mathcal{R}_3u(\log u)^{\beta}\right)(\nabla u)$. That is
\begin{align}\label{equ:8}
d\left(\mathcal{R}_1u+\mathcal{R}_2u^{\alpha}+\mathcal{R}_3u(\log u)^{\beta}\right)(\nabla u)&=\left(\mathcal{R}_1+\alpha\mathcal{R}_2u^{\alpha-1}+\mathcal{R}_3(\log u)^{\beta}+\beta\mathcal{R}_3(\log u)^{\beta-1}\right)F^2(\nabla u)\notag\\
&\quad+ud\mathcal{R}_1(\nabla u)+u^{\alpha}d\mathcal{R}_2(\nabla u)+u(\log u)^{\beta}d\mathcal{R}_3(\nabla u)\notag\\
&\le\left[ \mathcal{R}_1+\alpha\mathcal{R}_2u^{\alpha-1}+\mathcal{R}_3\left((\log u)^{\beta}+\beta(\log u)^{\beta-1}\right) \right]F^2(\nabla u)\notag\\
&\quad+\big[ uF_{\nabla u}(\nabla^{\nabla u}\mathcal{R}_1)+u^{\alpha}F_{\nabla u}(\nabla^{\nabla u}\mathcal{R}_2)\notag\\
&\quad+u(\log u)^{\beta}F_{\nabla u}(\nabla^{\nabla u}\mathcal{R}_3) \big]F(\nabla u)
\end{align}
Substituting \eqref{equ:8} into \eqref{equ:7}, for any $\phi\in C_{c}^{\infty}(M)$, we obtain
\begin{align}
&-\int_o^T\int_{M}2F(\nabla u)d\phi(\nabla^{\nabla u}(F(\nabla u)))d\mu dt-\int_0^T\int_{M}2\phi F(\nabla u)\partial_t(F(\nabla u))d\mu dt\notag\\
&\ge\int_0^T\int_{M}\phi\bigg\{-2(K+L_1)F^2(\nabla u)-2\left[ \mathcal{R}_1+\alpha\mathcal{R}_2u^{\alpha-1}+\mathcal{R}_3\left((\log u)^{\beta}+\beta(\log u)^{\beta-1}\right) \right]F^2(\nabla u)\notag\\
&\quad-2\big[ uF_{\nabla u}(\nabla^{\nabla u}\mathcal{R}_1)+u^{\alpha}F_{\nabla u}(\nabla^{\nabla u}\mathcal{R}_2)+u(\log u)^{\beta}F_{\nabla u}(\nabla^{\nabla u}\mathcal{R}_3)\big]F(\nabla u)\bigg\}d\mu dt,
\end{align}
which is equal to 
\begin{align}\label{equ:9}
&-\int_o^T\int_{M}d\phi(\nabla^{\nabla u}(F(\nabla u)))d\mu dt-\int_0^T\int_{M}\phi \partial_t(F(\nabla u))d\mu dt\notag\\
&\ge\int_0^T\int_{M}\phi\bigg\{-(K+L_1)F(\nabla u)-\left[ \mathcal{R}_1+\alpha\mathcal{R}_2u^{\alpha-1}+\mathcal{R}_3\left((\log u)^{\beta}+\beta(\log u)^{\beta-1}\right) \right]F(\nabla u)\notag\\
&\quad-\big[ uF_{\nabla u}(\nabla^{\nabla u}\mathcal{R}_1)+u^{\alpha}F_{\nabla u}(\nabla^{\nabla u}\mathcal{R}_2)+u(\log u)^{\beta}F_{\nabla u}(\nabla^{\nabla u}\mathcal{R}_3)\big] \bigg\}d\mu dt,
\end{align}
since $F(\nabla u)\geq 0$.
 On the other hand, we calculate $\mathcal{L}(u^2)$ and combine it with \eqref{1} and \eqref{equ:9} to obtain
\begin{align}\label{3.12}
\int_0^T\int_{M}\phi\mathcal{L}Hd\mu dt
&=\int_0^T\int_{M}\phi\bigg\{\mathcal{L}(F(\nabla u))+2F^2(\nabla u)+2u\mathcal{L}u\bigg\}d\mu dt\notag\\
&=\int_0^T\int_{M}\phi\bigg\{2F^2(\nabla u)-2\mathcal{R}_1u^2-2\mathcal{R}_2u^{\alpha+1}-2\mathcal{R}_3u^2(\log u)^{\beta}\bigg\}d\mu dt\notag\\
&\quad-\int_o^T\int_{M}d\phi(\nabla^{\nabla u}(F(\nabla u)))d\mu dt-\int_0^T\int_{M}\phi \partial_t(F(\nabla u))d\mu dt\notag\\
&\ge\int_0^T\int_{M}\phi\bigg\{-\left[K+L_1+ \mathcal{R}_1+\alpha\mathcal{R}_2u^{\alpha-1}+\mathcal{R}_3\left((\log u)^{\beta}+\beta(\log u)^{\beta-1}\right) \right]F(\nabla u)\notag\\
&\quad-[ uF_{\nabla u}(\nabla^{\nabla u}\mathcal{R}_1)+u^{\alpha}F_{\nabla u}(\nabla^{\nabla u}\mathcal{R}_2)+u(\log u)^{\beta}F_{\nabla u}(\nabla^{\nabla u}\mathcal{R}_3) ]\notag\\
&\quad+2F^2(\nabla u)-2\mathcal{R}_1u^2-2\mathcal{R}_2u^{\alpha+1}-2\mathcal{R}_3u^2(\log u)^{\beta}\bigg\}d\mu dt\notag\\
&\ge \int_0^T\int_{M}\phi\bigg\{(F(\nabla u)+u^2)^2-u^4\notag\\
&\quad-\left[2u^2+K+L_1+ \mathcal{R}_1+\alpha\mathcal{R}_2u^{\alpha-1}+\mathcal{R}_3\left((\log u)^{\beta}+\beta(\log u)^{\beta-1}\right) \right]F(\nabla u)\notag\\
&\quad-[ uF_{\nabla u}(\nabla^{\nabla u}\mathcal{R}_1)+u^{\alpha}F_{\nabla u}(\nabla^{\nabla u}\mathcal{R}_2)+u(\log u)^{\beta}F_{\nabla u}(\nabla^{\nabla u}\mathcal{R}_3) ]\notag\\
&\quad+F^2(\nabla u)-2\mathcal{R}_1u^2-2\mathcal{R}_2u^{\alpha+1}-2\mathcal{R}_3u^2(\log u)^{\beta}\bigg\}d\mu dt.
\end{align}
Using the assumptions and the inequality $xy\le \frac{1}{4}x^2+y^2$, we have
\begin{align}
&\left[2u^2+K+L_1+ \mathcal{R}_1+\alpha\mathcal{R}_2u^{\alpha-1}+\mathcal{R}_3\left((\log u)^{\beta}+\beta(\log u)^{\beta-1}\right) \right]F(\nabla u)\notag\\
&\le [2B^2+K^{+}+L_1^{+}+ \sigma_1+\vert\alpha\vert\sigma_1\max\{\delta^{\alpha-1},B^{\alpha-1}\}+\sigma_1(\tilde{D}_0+\beta\tilde{D}_1)]F(\nabla u) \notag\\
&\le \frac{1}{4}M_0^2+F^2(\nabla u),\label{equ:10}
\end{align}
where 
\begin{align*}
 	M_0&=2B^2+K^{+}+L_1^{+}+ \sigma_1+\vert\alpha\vert\sigma_1\hat{D}_{0}+\sigma_1(\tilde{D}_0+\beta\tilde{D}_1).\\
	\tilde{D}_i&=\max\left\{  \vert(\log \delta )^{\beta-i}\vert,\vert(\log B)^{\beta-i}\vert\right\},\; for \;\; i=0,1,  \\
	\hat{D}_j&=\max\left\{  \delta^{\alpha-1+j},B^{\alpha-1+j}\right\}, \; for \;\; j=0,1,2. 
\end{align*}
Moreover,
\begin{align}\label{equ:11}
&u^4+ uF_{\nabla u}(\nabla^{\nabla u}\mathcal{R}_1)+u^{\alpha}F_{\nabla u}(\nabla^{\nabla u}\mathcal{R}_2)+u(\log u)^{\beta}F_{\nabla u}(\nabla^{\nabla u}\mathcal{R}_3)\notag\\
&+2\mathcal{R}_1u^2+2\mathcal{R}_2u^{\alpha+1}+2\mathcal{R}_3u^2(\log u)^{\beta}\notag\\
\le& B^4+B\sigma_2+\hat{D}_1\sigma_2+B\tilde{D}_0\sigma_2+2\sigma_1B^2+2\sigma_1\hat{D}_2+2\sigma_1B^2\tilde{D}_0\notag\\
=:&B^4+M_1,
\end{align}
where $M_1=B\sigma_2+\hat{D}_1\sigma_2+B\tilde{D}_0\sigma_2+2\sigma_1B^2+2\sigma_1\hat{D}_2+2\sigma_1B^2\tilde{D}_0.$

Substituting \eqref{equ:10} and \eqref{equ:11} into \eqref{equ:9}, it implies that
\begin{align}\label{equ:12}
\mathcal{L}H\ge H^2-B^4-\frac{1}{4}M_0^2-M_1.
\end{align}
Fix an arbitrary time $t\in(0,T]$. Let $H$ attain its maximum at $(x_0,t_0)\in M\times [0,t]$. Since $H(x_0,t_0)>0$ (othertwise, the result is trivial), we know $t_0\in (0,t]$ and $H$ is positive on a neighborhood $U$ of $(x_0,t_0)$.
		Then, the RHS of \eqref{equ:12} must be nonpositive on $U$. Suppose, for the sake of contradiction, that it is positive on $U$, then $H$ is a local subsolution of 
		\begin{equation*}
			\mathcal{L}H\ge 0,
		\end{equation*}
		and its maximum should be attained at the boundary of $U$. This contradicts the fact that $(x_0,t_0)$ is an interior point of $U$. Thus, at $(x_0,t_0)$, we have
\begin{align}
0\ge\mathcal{L}H\ge  H^2-B^4-\frac{1}{4}M_0^2-M_1.
\end{align}
This reduces to solving a quadratic inequality in $H$. We find that $H\le \sqrt{B^4+\mathcal{G}}$, where $\mathcal{G}$ is defined as $\frac{1}{4}M_0^2+M_1$.
Furthermore, we obtain $$F(\nabla u)\le \sqrt{B^4+\mathcal{G}}-\delta^2\le B^2-\delta^2+\mathcal{G}^{\frac{1}{2}},$$
which completes the proof of the Shi-type gradient estimate. 
	\end{proof}
Next, we prove the Harnack inequality, which is a direct corollary of the Shi-type gradient estimates.
\begin{corollary}\label{Shi-Harnack}
Under the same assumptions and notation as in Theorem \ref{Shi-type}, it satisfies for any $x_1,x_2\in M$ and for any fixed $t_0\in(0,T)$ that
\begin{equation}
u(x_1, t_0) \le u(x_2, t_0) \cdot e^{C_3 d_F(x_1, x_2)},
\end{equation}
where $C_3$=$\frac{B^2-\delta^2+\mathcal{G}^{\frac{1}{2}}}{\delta}$. The Finsler distance is defined by $ d_F(x_1, x_2)=\int_0^1 F(\dot{\gamma}(s))ds$, where $\gamma(s)$ is a  minimal geodesic from $x_2$ to $x_1$ in $B_p(R)$.
\end{corollary}

\begin{proof}[Proof of Corollary]\eqref{Shi-Harnack}
Based on the result of Theorem \eqref{Shi-type}, we have 
$$
F(\nabla \log u)\le \frac{B^2-\delta^2+\mathcal{G}^{\frac{1}{2}}}{\delta}.
$$ 
Let $\gamma(s)$ be a  minimal geodesic in $B_p(R)$ parameterized by arc-length, connecting from the point $x_2$ to $x_1$. According to the definition of the Finsler distance $ d_F(x_1, x_2)=\int_0^1 F(\dot{\gamma}(s))ds $, by setting $f(s) = \log u(\gamma(s))$, we have
\begin{align*}
f(x_1,t_0) - f(x_2,t_0) &= \int_0^1 \frac{d}{ds}[f(\gamma(s))] ds \\
&= \int_0^1 df(\dot{\gamma}(s)) ds  \\
&\le \int_0^1F(\nabla f)F(\dot{\gamma}(s)) ds  \\
&\le C_3\int_0^1 F(\dot{\gamma}(s))ds\\
&=C_3d_F(x_1,x_2),
\end{align*}
where $C_3=\frac{B^2-\delta^2+\mathcal{G}^{\frac{1}{2}}}{\delta}$.
Exponentiating both sides of the inequality (since the exponential function is increasing), we conclude that
$$
u(x_1, t_0) \le u(x_2, t_0) \cdot e^{C_3d_F(x_1, x_2)}.
$$
\end{proof}
	\section{Hamilton-type gradient estimate}\label{sec4}
In this section, we further establish the Hamilton-type gradient estimates and the Harnack inequality.

\begin{proof}[Proof of Theorem \ref{Hamilton-type}]
First, set $H=\sqrt{1+\log \frac{B}{u}}$=$\sqrt{\log \frac{D}{u}}\ge1$, where $D=eB$.
Consequently, $u=De^{-H^2}$,  $u_t=D(-2H_t)e^{-H^2}H$, $\nabla u=D(-2H)e^{-H^2}\nabla^{\nabla u}H,$ and
$$ \Delta u=-2DHe^{-H^2}\Delta^{\nabla u}H-2DF_{\nabla u}^2(\nabla^{\nabla u}H)e^{-H^2}+4DH^2e^{-H^2}F_{\nabla u}^2(\nabla^{\nabla u}H).$$
According to the definition of the operator $\mathcal{L}$ (cf. \eqref{L}), \eqref{2} becomes
\begin{equation}
\mathcal{L}H=F_{\nabla u}^2(\nabla^{\nabla u}H)(2H-\frac{1}{H})+\frac{\mathcal{R}_1}{2H}+\frac{\mathcal{R}_2}{2H}(De^{-H^2})^{\alpha-1}+\frac{\mathcal{R}_3}{2H}(\log D-H^2)^{\beta}.
\end{equation}

Now, let $\omega=F^2(\nabla^{\nabla u}H)$. By applying the Bochner formula \eqref{lem2.1} and Lemma \ref{lem1}, we obtain 
\begin{align}\label{equ:15}
\mathcal{L}\omega&=\Delta^{\nabla u}F^2(\nabla^{\nabla u}H)-\partial_tF^2(\nabla^{\nabla u}H)\notag\\
&=2\big(dH(\nabla^{\nabla u}(\Delta^{\nabla u} H))+Ric^{\infty}(\nabla^{\nabla u} H)+\Vert\nabla^2H\Vert_{HS(\nabla u)}^2\big)\notag\\
&\quad-2h(\nabla^{\nabla u} H)(\nabla^{\nabla u} H,\nabla^{\nabla u} H)-2dH_t(\nabla^{\nabla u} H).\notag\\
&=2Ric^{\infty}(\nabla^{\nabla u}H)+2\Vert\nabla^2H\Vert_{HS(\nabla u)}^2-2h(\nabla^{\nabla u} H)(\nabla^{\nabla u} H,\nabla^{\nabla u} H)\notag\\
&\quad+2\big(dH(\nabla^{\nabla u}(\Delta^{\nabla u} H))-dH_t(\nabla^{\nabla u} H)\big).
\end{align}
 
Direct  computation yields that
\begin{align}\label{equ:16}
&dH\big(\nabla^{\nabla u}(\Delta^{\nabla u} H)\big)-dH_t\big(\nabla^{\nabla u} H\big)\notag\\
=&dH\bigg(\nabla^{\nabla u}\big( \partial_tH-F_{\nabla u}^2(\nabla^{\nabla u}H)(\frac{1}{H}-2H)+\frac{\mathcal{R}_1}{2H}\notag\\
&+\frac{\mathcal{R}_2}{2H}(De^{-H^2})^{\alpha-1}+\frac{\mathcal{R}_3}{2H}(\log D-H^2)^{\beta}\big)\bigg)-dH_t(\nabla^{\nabla u}H)\notag\\
=&-dH\Bigg(\nabla^{\nabla u}\Big(F_{\nabla u}^2(\nabla^{\nabla u}H)(\frac{1}{H}-2H)\Big)\Bigg)+dH\left(\nabla^{\nabla u}\big(\frac{\mathcal{R}_1}{2H}\big)\right)\notag\\
&\quad+dH\Bigg(\nabla^{\nabla u}\Big(\frac{\mathcal{R}_2}{2H}(De^{-H^2})^{\alpha-1}\Big)\Bigg)+dH\Bigg(\nabla^{\nabla u}\Big(\frac{\mathcal{R}_3}{2H}(\log D-H^2)^{\beta}\Big)\Bigg),
\end{align}
where we can further deduce that
\begin{align}
dH\Bigg(\nabla^{\nabla u}\Big(F_{\nabla u}^2(\nabla^{\nabla u}H)(\frac{1}{H}-2H)\Big)\Bigg)&=(\frac{1}{H}-2H)dH(\nabla^{\nabla u}F_{\nabla u}^2(\nabla^{\nabla u}H))\notag\\
&\quad-(\frac{1}{H^2}+2)F_{\nabla u}^4(\nabla^{\nabla u}H),\label{equ:17}\\
dH\Bigg(\nabla^{\nabla u}\big(\frac{\mathcal{R}_1}{2H}\big)\Bigg)&=\frac{1}{2H}dH(\nabla^{\nabla u}\mathcal{R}_1)-\frac{\mathcal{R}_1}{2H^2}F_{\nabla u}^2(\nabla^{\nabla u}H),\label{equ:18}\\
dH\Bigg(\nabla^{\nabla u}\Big(\frac{\mathcal{R}_2}{2H}(De^{-H^2})^{\alpha-1}\Big)\Bigg)&=\frac{(De^{-H^2})^{\alpha-1}}{2H}dH(\nabla^{\nabla u}\mathcal{R}_2)\notag\\
&\quad-\frac{(De^{-H^2})^{\alpha-1}}{2H^2}\mathcal{R}_2F_{\nabla u}(\nabla^{\nabla u}H)\notag\\
&\quad+(1-\alpha)\mathcal{R}_2(De^{-H^2})^{\alpha-1}F_{\nabla u}^2(\nabla^{\nabla u}H),\label{equ:19}
\end{align}
and
\begin{align}\label{equ:20}
dH\Bigg(\nabla^{\nabla u}\Big(\frac{\mathcal{R}_3}{2H}(\log D-H^2)^{\beta}\Big)\Bigg)&=\frac{1}{2H}(\log D-H^2)^{\beta}dH(\nabla^{\nabla u}\mathcal{R}_3)\notag\\
&\quad-\frac{\mathcal{R}_3}{2H^2}(\log D-H^2)^{\beta}F_{\nabla u}^2(\nabla^{\nabla u}H)\notag\\
&\quad-\beta(\log D-H^2)^{\beta-1}\mathcal{R}_3F_{\nabla u}^2(\nabla^{\nabla u}H).
\end{align}
By substituting \eqref{equ:16}-\eqref{equ:20} with \eqref{equ:15}, we obtain on $M_u$ that
\begin{align}
\mathcal{L}\omega&=2Ric^{\infty}(\nabla^{\nabla u}H)+2\Vert\nabla^2H\Vert_{HS(\nabla u)}^2)-2h_{\nabla H}(\nabla^{\nabla u} H,\nabla^{\nabla u} H)\notag\\
&\quad-2(\frac{1}{H}-2H)dH(\nabla^{\nabla u}F_{\nabla u}^2(\nabla^{\nabla u}H))+2(\frac{1}{H^2}+2)F_{\nabla u}^4(\nabla^{\nabla u}H)\notag\\
&\quad+\frac{1}{H}dH(\nabla^{\nabla u}\mathcal{R}_1)-\frac{\mathcal{R}_1}{H^2}F_{\nabla u}^2(\nabla^{\nabla u}H)\notag\\
&\quad+\frac{(De^{-H^2})^{\alpha-1}}{H}dH(\nabla^{\nabla u}\mathcal{R}_2)-\frac{(De^{-H^2})^{\alpha-1}}{H^2}\mathcal{R}_2F_{\nabla u}(\nabla^{\nabla u}H)\notag\\
&\quad+2(1-\alpha)\mathcal{R}_2(De^{-H^2})^{\alpha-1}F_{\nabla u}^2(\nabla^{\nabla u}H)+\frac{1}{H}(\log D-H^2)^{\beta}dH(\nabla^{\nabla u}\mathcal{R}_3)\notag\\
&\quad-\frac{\mathcal{R}_3}{H^2}(\log D-H^2)^{\beta}F_{\nabla u}^2(\nabla^{\nabla u}H)-2\beta(\log D-H^2)^{\beta-1}\mathcal{R}_3F_{\nabla u}^2(\nabla^{\nabla u}H)\notag\\
&\le -2\rho^2(K+L_1)\omega+\Bigl(2(1-\alpha)\mathcal{R}_2(De^{-H^2})^{\alpha-1}-2\beta(\log D-H^2)^{\beta-1}\mathcal{R}_3\Bigr)\omega\notag\\
&\quad+\frac{1}{H}\Bigl(dH(\nabla^{\nabla u}\mathcal{R}_1)+(De^{-H^2})^{\alpha-1}dH(\nabla^{\nabla u}\mathcal{R}_2)+(\log D-H^2)^{\beta}dH(\nabla^{\nabla u}\mathcal{R}_3)\Bigr)\notag\\
&\quad-\frac{1}{H^2}\Bigl( \mathcal{R}_1+\mathcal{R}_2(De^{-H^2})^{\alpha-1}+\mathcal{R}_3(\log D-H^2)^{\beta} \Bigr)\omega\notag\\
&\quad+2(2+\frac{1}{H^2})\omega^2+2(2H-\frac{1}{H})dH(\nabla^{\nabla u}\omega).
\end{align}


		Thus, the RHS of \eqref{equ:12} is nonpositive on $U$. If not, $\omega$ is a local weak subsolution to 
		\begin{equation*}
			\mathcal{L}\omega\ge 0,
		\end{equation*}
		and its maximum must be attained at the boundary of $U$. This contradicts the fact that $(x_0,t_0)$ is an interior point of $U$. Therefore, at $(x_0,t_0)$, we have that 
\begin{align}
0\ge\mathcal{L}\omega&\ge- 2\rho^2(K+L_1) \omega+\Bigl(2(\alpha-1)\mathcal{R}_2(De^{-H^2})^{\alpha-1}+2\beta(\log D-H^2)^{\beta-1}\mathcal{R}_3\Bigr)\omega\notag\\
&\quad+\frac{1}{H^2}\Bigl( \mathcal{R}_1+\mathcal{R}_2(De^{-H^2})^{\alpha-1}+\mathcal{R}_3(\log D-H^2)^{\beta} \Bigr)\omega\notag\\
&\quad+\frac{1}{H}\Bigl(dH(\nabla^{\nabla u}\mathcal{R}_1)+(De^{-H^2})^{\alpha-1}dH(\nabla^{\nabla u}\mathcal{R}_2)+(\log D-H^2)^{\beta}dH(\nabla^{\nabla u}\mathcal{R}_3)\Bigr)\notag\\
&\quad+2\left(\frac{2H^2+1}{H^2}\right)\omega^2.
\end{align}
Multiplying both sides by $-\frac{H^2}{2H^2+1}$ shows that
\begin{align}\label{*}
2\omega^2 &\le \frac{2 H^2}{2H^2+1}(\rho^2K+L_1)\omega+\frac{2H^2}{2H^2+1}(\alpha-1)\mathcal{R}_2(De^{-H^2})^{\alpha-1}\omega\notag\\
&\quad+\frac{2 H^2}{2H^2+1}\beta(\log D-H^2)^{\beta-1}\mathcal{R}_3\omega\notag\\
&\quad+\frac{1}{2H^2+1}\Bigl(\mathcal{R}_1+\mathcal{R}_2(De^{-H^2})^{\alpha-1}+\mathcal{R}_3(\log D-H^2)^{\beta}\Bigr)\omega\notag\\
&\quad-\frac{H}{2H^2+1}\Bigl(dH(\nabla^{\nabla u}\mathcal{R}_1)+(De^{-H^2})^{\alpha-1}dH(\nabla^{\nabla u}\mathcal{R}_2)+(\log D-H^2)^{\beta}dH(\nabla^{\nabla u}\mathcal{R}_3)\Bigr).
\end{align}

It follows from the inequalities $0<\frac{2H^2}{2H^2+1}\le 1$, $0<\frac{1}{2H^2+1}\le \frac{1}{3}$ and $\frac{H}{1+2H^2}\le \frac{1}{3}$ that
\begin{align}
\frac{2 H^2}{2H^2+1}(\alpha-1)\mathcal{R}_2(De^{-H^2})^{\alpha-1} &\le \vert\alpha-1\vert\sigma_1\hat{D}_0,\label{equ:21}\\
\frac{2 H^2}{2H^2+1}\beta(\log D-H^2)^{\beta-1}\mathcal{R}_3S&\le \beta\sigma_1 \tilde{D}_1,\label{equ:22}
\end{align}
and
\begin{align}\label{equ:23}
\frac{1}{2H^2+1}\Bigl(\mathcal{R}_1+\mathcal{R}_2(De^{-H^2})^{\alpha-1}+\mathcal{R}_3(\log D-H^2)^{\beta}\Bigr) \le \frac{\sigma_1}{3}(1+\hat{D}_0+\tilde{D}_0).
\end{align}

By the Cauchy-Schwarz inequality, we deduce that 
\begin{align}\label{equ:24}
&-\frac{ H}{2H^2+1}\Bigl(dH(\nabla^{\nabla u}\mathcal{R}_1)+(De^{-H^2})^{\alpha-1}dH(\nabla^{\nabla u}\mathcal{R}_2)+(\log D-H^2)^{\beta}dH(\nabla^{\nabla u}\mathcal{R}_3)\Bigr)\notag\\
\le& \frac{1}{3}\Bigl(\epsilon_1\omega^2+\frac{1}{\epsilon_1}F_{\nabla u}(\nabla^{\nabla u}\mathcal{R}_1)+\epsilon_2\omega^2+\frac{1}{\epsilon_2}\sup\{u^{\alpha-1} \}F_{\nabla u}(\nabla^{\nabla u}\mathcal{R}_2) \notag\\
&\,\,\quad+\epsilon_3\omega^2+\frac{1}{\epsilon_3}\sup\{ \vert \log u \vert^{\beta}\}F_{\nabla u}(\nabla^{\nabla u}\mathcal{R}_3)\Bigr)\notag\\
\le& \frac{1}{3}(\epsilon_1+\epsilon_2+\epsilon_3)\omega^2+\frac{\sigma_2}{3}(\frac{1}{\epsilon_1}+\frac{1}{\epsilon_2}\hat{D}_0+\frac{1}{\epsilon_3}\tilde{D}_0),
\end{align}
where $\epsilon_i (i=1,2,3)$ are positive constants to be determined later.
Substituting \eqref{equ:21}-\eqref{equ:24} into \eqref{*} gives that
\begin{align}\label{4.16}
    \Bigl(2-\frac{1}{3}(\epsilon_1+\epsilon_2+\epsilon_3)\Bigr)\omega^2 &\le \Biggl\{\rho^2(K+L_1)+ \Bigl(\vert\alpha-1\vert\sigma_1\hat{D}_0+\beta\sigma_1 \tilde{D}_1+\frac{\sigma_1}{3}(1+\hat{D}_0+\tilde{D}_0)\Bigr)\Biggr\}\omega\notag\\
&\quad+\frac{\sigma_2}{3}(\frac{1}{\epsilon_1}+\frac{1}{\epsilon_2}\hat{D}_0+\frac{1}{\epsilon_3}\tilde{D}_0).
\end{align}
Recall that $a\le a_2+\sqrt{a_1}$ holds, provided $a^2\le a_1+a_2a$ with $a,a_1,a_2>0$. By choosing $\epsilon_1+\epsilon_2+\epsilon_3=3$, \eqref{4.16} implies that
\begin{equation}
    \omega\le C_1+C_2,
\end{equation}
where
\begin{align*}
    C_1&=\frac{\sigma_2}{3}(\frac{1}{\epsilon_1}+\frac{1}{\epsilon_2}\hat{D}_0+\frac{1}{\epsilon_3}\tilde{D}_0),\\
    C_2&=\sqrt{\rho^2(K+L_1)+ \left[\vert\alpha-1\vert\sigma_1\hat{D}_0+\beta\sigma_1 \tilde{D}_1+\frac{\sigma_1}{3}(1+\hat{D}_0+\tilde{D}_0)\right]}.
\end{align*}
Hence,
\begin{equation}\label{eq-4.17}
       F(-\nabla \log u)\le 2\sqrt{(1+\log\frac{B}{\delta})(C_1+C_2)},
\end{equation}
and more,
\begin{equation}
  F(\nabla \log u)\le 2\rho\sqrt{(1+\log\frac{B}{\delta})(C_1+C_2)}.
\end{equation}
where $\rho$ is given by \eqref{rho}. Furthermore,
 \begin{equation}
        \max\{F(\nabla \log u), F(\nabla -\log u)\}\le 2\rho\sqrt{(1+\log\frac{B}{\delta})(C_1+C_2)},
    \end{equation}
This completes the proof of Theorem \ref{Hamilton-type}.
\end{proof}
Besides, we prove the Harnack inequality, which is a direct corollary of the Hamilton-type gradient estimates.
\begin{corollary}\label{Hamilton-Harnack}
Under the same assumptions and notation as in Theorem \ref{Shi-type}, it satisfies for any $x_1,x_2\in M$ and for any fixed $t_0\in(0,T)$ that
\begin{equation}
u(x_1, t_0) \le u(x_2, t_0) \cdot e^{C_4 d_F(x_1, x_2)},
\end{equation}
where $C_4=2\sqrt{(1+\log\frac{B}{\delta})(C_1+C_2)},$ with $C_1,C_2$ given in Theorem \ref{Hamilton-type}.  The Finsler distance is defined by $ d_F(x_1, x_2)=\int_0^1 F(\dot{\gamma}(s))ds$, where $\gamma(s)$ is a  minimal geodesic from $x_2$ to $x_1$ in $B_p(R)$.
\end{corollary}
\begin{proof}[Proof of Corollary \ref{Hamilton-Harnack}] 
Following the same approach as in the proof of \eqref{Hamilton-Harnack}, and according to \eqref{eq-4.17}, we have
\begin{align*}
f(x_1,t_0) - f(x_2,t_0) &= \int_0^1 \frac{d}{ds}[f(\gamma(s))] ds \\
&= \int_0^1 df(\dot{\gamma}(s)) ds  \\
&\le \int_0^1F(-\nabla f)F(\dot{\gamma}(s)) ds  \\
&\le C_4\int_0^1 F(\dot{\gamma}(s))ds\\
&=C_4d_F(x_1,x_2),
\end{align*}
where $C_4=2\sqrt{(1+\log\frac{B}{\delta})(C_1+C_2)}$.
Exponentiating both sides of the inequality (since the exponential function is increasing), we conclude that
$$
u(x_1, t_0) \le u(x_2, t_0) \cdot e^{C_4d_F(x_1, x_2)}.
$$
\end{proof}

\section*{Acknowledgments}

The authors are grateful to the reviewers for their careful reviews and valuable comments. 

\section*{Declarations}

\subsection*{Ethical Approval}
Ethical Approval is not applicable to this article as no human or animal studies in this study.

\subsection*{Funding} 
The second author is supported partially by the NNSFC (Nos. 12001099, 12271093).

\subsection*{Data availability statement}
Data sharing is not applicable to this article as no new data were created or analyzed in this study.

\subsection*{Materials availability statement}
Materials sharing is not applicable to this article.

\subsection*{Conflict of interest/Competing interests}
All authors disclosed no relevant relationships.
\end{document}